\documentclass{article}[10pt]

\usepackage{amssymb}

\usepackage[english,francais]{babel}

\newtheorem{e-proposition}[theorem]{Proposition}

\newtheorem{e-definition}[theorem]{Definition\rm}

\newtheorem{theoreme}{Th\'eor\`eme}[section]

\newtheorem{definition}[theoreme]{D\'efinition\rm}

\def\og{\leavevmode\raise.3ex\hbox{$\scriptscriptstyle\langle\!\langle$~}}
\def\fg{\leavevmode\raise.3ex\hbox{~$\!\scriptscriptstyle\,\rangle\!\rangle$}}

\begin{document}
\selectlanguage{francais}

\noindent\textbf{\large{Produit de Moyal stochastique sur l'espace de Wiener}}

\vskip10mm

\noindent Giuseppe Dito \& R\'emi L\'eandre

\vskip5mm

\noindent{Institut de Math\'ematiques de Bourgogne\\CNRS UMR 5584\\ 
Universit\'e de Bourgogne\\ B.P. 47870\\ 21078 Dijon CEDEX, France
\vskip10mm
\noindent e-mail : \texttt{giuseppe.dito@u-bourgogne.fr, remi.leandre@u-bourgogne.fr}
}
\vskip10mm
\noindent{\bf Abstract}
\vskip 0.5\baselineskip
\noindent
We propose a stochastic extension of deformation quantization on a Hilbert space.
The Moyal product is defined in this context on the space of functionals belonging
to all of the Sobolev spaces of the Malliavin calculus.
\vskip10mm
\selectlanguage{francais}
\noindent{\bf R\'esum\'e}
\vskip 0.5\baselineskip
\noindent
Nous proposons une version stochastique de la quantification par déformation sur un espace de Hilbert.
Dans ce cadre, le produit de Moyal est d\'efini sur l'espace des fonctionnelles appartenant \`a tous les
espaces de Sobolev du Calcul de Malliavin.


\selectlanguage{francais}

\section{Introduction}\label{sect1}

Le premier de cette note auteur a r\'ecemment \'etudi\'e la th\'eorie de la quantification par d\'eformation 
sur un espace de Hilbert \cite{Di} qui est une
extension en dimension infinie de \cite{BFFLS} (Voir \cite{DiS} pour un article de revue r\'ecent). 
Les formules obtenues sont tr\`es analogues \`a celles de l'analyse de Wiener sur un espace de Wiener abstrait 
$H \subset W$ \cite{IW}, \cite{Ma}, \cite{Nu}, \cite{Us}. On remarque n\'eanmoins que les sp\'ecificit\'es du Calcul de Malliavin 
par rapport aux travaux ant\'erieurs (Voir les travaux de Fomin, Albeverio, Elworthy, Hida, Berezanskii) sont les suivantes :

a) il consid\`ere l'intersection de tous les espaces $L^p$ sur l'espace de Wiener au lieu de consid\'erer uniquement l'espace 
de Hilbert $L^2$, si bien que l'ensemble des fonctions test du Calcul de Malliavin constitue une {\bf{alg\`ebre}}
gr\^ace \`a l'in\'egalit\'e de H\"older ;

b) il effectue la compl\'etion des op\'erations diff\'erentielles connues en ce temps sur l'espace de Wiener, 
si bien qu'une fonctionnelle appartenant \`a tous les espaces de Sobolev du Calcul de Malliavin n'est en g\'en\'eral 
que presque surement d\'efinie, puisqu'il n'y a pas de th\'eor\`eme d'injection de Sobolev en dimension infinie.

Par rapport au travail \cite{Di}, nous consid\'erons l'alg\`ebre des fonctionnelles test au sens de Malliavin. 
Les formules de \cite{Di} restent valides, mais seulement presque surement, et avec des conditions naturelles d'int\'egrabilit\'e.

Comme l'a remarqu\'e Meyer \cite{Me}, on peut se ramener au cas de l'\'etude de $L^2([0,1])$ et de l'espace de Wiener usuel.

\section{Brefs rappels sur le calcul de Malliavin}\label{sect2}

Soit $H \subset W$ un espace de Wiener abstrait \cite{Nu}. Soit $F$ une fonctionnelle appartenant \`a tous les 
espaces de Sobolev du Calcul de Malliavin. $\nabla^rF$ r\'ealise un \'el\'ement al\'eatoire de $H^{\hat{\otimes}r}$, 
$\hat{\otimes}$ d\'esignant le produit tensoriel hilbertien sym\'etrique, 
puisqu'on ne d\'erive que dans la direction de $H$. On pose
\begin{equation}\label{eq21}
\Vert F \Vert_{r,p} = E[\Vert \nabla^rF\Vert^p]^{1/p}.
\end{equation}
$W_{r,p}$ est l'espace des fonctionnelles telles que $\Vert F\Vert_{r,p} < \infty$. $W_{\infty-}$ est 
l'espace des fonctionnelles telles que $\Vert F \Vert_{r,p} < \infty$ pour tout $r$ et tout $p$. 
C'est une {\bf{alg\`ebre commutative}}. De fa\c con analogue, on peut d\'efinir des espaces de 
Sobolev de fonctionnelles \`a valeurs dans des espaces de Hilbert.

Dans le cas o\`u $H = L^2([0,1]) \subset C([0,1])$ est l'espace de Wiener habituel
 (on assimile $h \in L^2([0,1])$  \`a la fonction $s \mapsto \int_0^sh(u) du$). $\nabla^rF$ 
est un \'el\'ement al\'eatoire $F(s_1,\ldots,s_r)$. La fonctionnelle $F$ appartient \`a $W_{r,p}$ si
\begin{equation}\label{eq22}
E[(\int_{[0,1]^r}\vert \nabla^rF(s_1,\ldots,s_r)\vert^2ds_1\cdots ds_r)^{p/2}]^{1/p} < \infty. 
\end{equation}

\section{Quantification par d\'eformation d'un espce de Wiener abstrait}\label{sect3}

Soit $H \subset W$ un espace de Wiener abstrait. Au lieu de consid\'erer comme dans \cite{Di} 
l'alg\`ebre des fonctions $C^\infty$ au sens de Fr\'echet sur $H$, on prend $W_{\infty-}$ comme 
alg\`ebre de fonctions d\'efinies sur $W$.

\begin{definition}\label{def31}
$(W_{\infty-}, \{\cdot,\cdot\})$ est appel\'e un espace de Wiener de Poisson si les conditions suivantes sont r\'ealis\'ees :

i) il existe une application bilin\'eaire al\'eatoire $P$ de $H\times H$ \`a valeurs dans $\mathbb{R}$ appartenant 
\`a tous les espaces de Sobolev sur $W$ ;

ii) si on pose pour $F$ et $G$ appartenant \`a tous les espaces de Sobolev du Calcul de Malliavin :

\begin{equation}\label{eq31}
\{F,G\} = P(\nabla F, \nabla G),
\end{equation}
alors $\{\cdot,\cdot\}$ est un crochet de Poisson.
\end{definition}

Un crochet de Poisson sur l'espace de Wiener est un cas particulier d'op\'erateur $r$-diff\'erentiel
sur l'espace de Wiener d\'efini par :

\begin{definition}\label{def32}
 Un op\'erateur $r$-diff\'erentiel $A$ sur $W_{\infty-}$ est d\'efini par les donn\'ees suivantes :

i) une famille finie $a^{n_1,\dots,n_r}$ d'applications de $H^{\hat{\otimes}n_1}\times \cdots\times H^{\hat{\otimes}n_r}$ 
\`a valeurs dans $\mathbb{R}$ appartenant \`a tous les espaces de Sobolev du Calcul de Malliavin;

ii) pour $F_1,\ldots,F_r \in W_{\infty-}$, on a :
\begin{equation}\label{eq32}
A(F_1,\ldots,F_r) = \sum a^{n_1,\ldots,n_r}(\nabla^{n_1}F_1,\ldots,\nabla^{n_r}F_r).
\end{equation}
\end{definition}
Il r\'esulte de l'in\'egalit\'e de H\"older que $A$ applique continuement $W_{\infty-}^r$ sur $W_{\infty-}$.

Nous pouvons maintenant donner la d\'efinition d'une quantification par d\'efor\-mation d'un espace de Wiener abstrait.

\begin{definition}\label{def33}
Soit $(W_{\infty-},\{\cdot,\cdot\})$ un espace de Wiener de Poisson. 
On consid\`ere l'espace $W_{\infty-}[[h]]$ des s\'eries formelles $\sum h^n F_n$ o\`u $F_n$ appartient \`a $W_{\infty-}$. 
Un star-produit  $\star_\hbar$ est la donn\'ee d'une application $\mathbb{R}[[\hbar]]$-bilin\'eaire de 
$W_{\infty-}[[\hbar]]\times W_{\infty-}[[\hbar]]$ dans $W_{\infty-}[[h]]$ telle que :

i) $F\star_\hbar G = \sum \hbar^rC_r(F,G)$ ;

ii) $C_0(F,G) = FG$ ;

iii) $C_1(F,G)-C_1(G,F) = 2\{F,G\}$ ;

iv) $C_r$ est un op\'erateur bidiff\'erentiel au sens de la d\'efinition~\ref{def32} ;

v) $F\star_\hbar(G\star_\hbar H) = (F\star_\hbar G)\star_\hbar H$.
\end{definition}

La partie suivante est consacr\'ee \`a l'\'etude d'un exemple.

\section{Le produit de Moyal stochastique}\label{sect4}

Nous consid\'erons $H = L^2([0,1]) \oplus L^2([0,1])$ et $W =C([0,1))\oplus C([0,1])$. 
Le symbole $w^1$ 
d\'esigne un \'el\'ement du premier espace de Wiener et $w^2$ un \'el\'ement du deuxieme  
(les deux espaces de Wiener dans $W$ sont ind\'ependants). 
Si $F$ et $G$ sont deux fonctionnelles appartenant \`a $W_{\infty-}$ 
pour l'espace de Wiener {\bf{total}}, on introduit :
\begin{equation}\label{eq41}
\{F,G\}(w^1,w^2) = \int_0^1\nabla_1F(s)\nabla_2G(s)ds-\int_0^1\nabla_1G(s)\nabla_2F(s)ds.
\end{equation}
Dans la premi\`ere int\'egrale, nous d\'erivons $F$ par rapport au premier espace et $G$ 
par rapport au second espace de Hilbert, et r\'eciproquement dans la seconde int\'egrale. 

Nous avons clairement :

\begin{theoreme} 
$W_{\infty-}$ muni du crochet (\ref{eq41}) est un espace de Wiener de Poisson.
\end{theoreme}

Soit $\alpha_1,\ldots,\alpha_r,\beta_1,\ldots,\beta_r$ des entiers \'egaux \`a 1 ou 2. 
Nous posons (Voir \cite{Di} (6)) :
\begin{eqnarray}\label{eq42}
& &<\nabla^r_{\alpha_1,\ldots,\alpha_r}F,\nabla^r_{\beta_1,\ldots,\beta_r}G> 
= \\
& &\qquad \int_{[0,1]^r}ds_1\cdots ds_r\nabla^r_{\alpha_1,\ldots,\alpha_r}F(s_1,\ldots,s_r)
\nabla^r_{\beta_1,\ldots,\beta_r}G(s_1,\ldots,s_r).\nonumber
\end{eqnarray}
Cela d\'efinit clairement un op\'erateur bidiff\'erentiel au sens de la d\'efinition~\ref{def32}.

Soit $\Lambda$ la  matrice symplectique $2\times 2$ d\'efinie par $\Lambda^{1,2} = 1$. Suivant \cite{Di} (7), nous posons :
\begin{equation}\label{eq43}
C_r(F,G) = \sum_{\alpha_1,\ldots,\alpha_r,\beta_1,\ldots,\beta_r}
\Lambda^{\alpha_1,\beta_1}\cdots\Lambda^{\alpha_r,\beta_r}<\nabla^r_{\alpha_1,\dots,\alpha_r}F,
\nabla^r_{\beta_1,\ldots,\beta_r}G>.
\end{equation}
Comme dans (\ref{eq43}) la somme est finie, $C_r$ est un op\'erateur bidiff\'erentiel sur l'espace de Wiener total.

\begin{definition}
 Le produit de Moyal stochastique est d\'efini par :
\begin{equation}\label{eq44}
F*_h G = FG + \sum_{r\geq 1}{h^r\over r!}C_r(F,G).
\end{equation}
\end{definition}
En proc\'edant comme dans \cite{Di}, Th\'eoreme 1, on montre :

\begin{theoreme} 
Le produit de Moyal stochastique r\'ealise une quantification par d\'eformation de 
l'espace de Wiener de Poisson $(W_{\infty-},\{\cdot,\cdot\})$.
\end{theoreme}

\section{L'espace de phase de l'espace de Wiener}\label{sect5}

On a proc\'ed\'e ci-dessus comme si $W$ \'etait l'espace de phase de 
l'espace de Wiener $C([0,1])$ muni de sa forme symplectique canonique. 
Mais comme l'ont montr\'e les travaux de L\'eandre \cite{L1,L2} 
il faut plut\^ot consid\'erer comme espace de phase de l'espace de Wiener l'espace 
$C([0,1]) \oplus L^2([0,1])$ o\`u on a assimil\'e $L^2([0,1])$ \`a l'espace de Sobolev 
$H^{1,2}([0,1])$ des applications $f$ telles $\int_0^1\vert d/dsf(s)\vert^2ds <\infty$. 
Sur $L^2([0,1])$ nous mettons la mesure Gaussienne avec espace auto-reproduisant $H^{1,2}$ 
des $h$ tels que $\int_0^1\vert d/ds h(s)\vert^2ds = \Vert h\Vert_{1,2}^2< \infty$. 
Si l'on prend la d\'eriv\'ee d'une fonctionnelle $F$ dans la direction de $H_{1,2}$,
nous assimilons son noyau al\'eatoire $\nabla_2F(s)$ \`a $\int_0^{s}\nabla_2F(t)dt$. 
En utilisant cette identification, nous pouvons donner une extension \`a ce cadre de la forme  
symplectique canonique sur l'espace de phase $W = C([0,1]) \oplus L^2([0,1])$ en utilisant 
la d\'erivation suivant $L^2([0,1])$ pour $C([0,1])$ et suivant $H^{1,2}$  pour $L^2([0,1])$.

La formule (\ref{eq42}) s'interpr\`ete facilement dans ce nouveau formalisme et permet 
d'obtenir un nouveau produit de Moyal stochastique.


\begin{thebibliography}{00}

\bibitem{BFFLS} F. Bayen, M. Flato, C. Fronsdal, A. Lichnerowicz, D. Sternheimer, 
                 Deformation theory and quantization, Ann. Phys. 111 (1978) 61--151.

\bibitem{Di} G. Dito, 
             Deformation quantization on a Hilbert space, Preprint: math.QA/0406583 (2004).

\bibitem{DiS} G. Dito, D. Sternheimer, 
              Deformation quantization: genesis, developments and metamorphoses, in: G. Halbout (ed.),
              Deformation quantization,  IRMA Lectures in Maths. Theor. Phys. Walter de Gruyter, 2002 pp. 9--54.

\bibitem{IW}  N. Ikeda, S. Watanabe, 
              Stochastic differential equations and diffusion processes, 2nd ed., North Holland, 1989.

\bibitem{L1} L\'eandre R.: 
             Stochastic gauge transform of the string bundle, J. Geom. Phys. 26 (1998) 1--25.

\bibitem{L2} L\'eandre R.: 
             Stochastic cohomology of the frame bundle of the loop space, J. Nonlinear. Math. Phys. 5 (1998) 17--31.

\bibitem{Ma} P. Malliavin, 
             Stochastic analysis, Springer-Verlag, 1997.

\bibitem{Me} P. A. Meyer,
             Quantum probability for probabilists, Lect. Notes Math. 1538, Springer-Verlag, 1993.

\bibitem{Nu} D. Nualart,
             Malliavin Calculus and related topics, Springer-Verlag, 1995.

\bibitem{Us} A. S. \"Ust\"unel,
             An introduction to analysis on Wiener space, Lect. Notes Math. 1610, Springer-Verlag, 1995.

\end{thebibliography}
\end{document}